\documentclass[preprint,12pt]{elsarticle}

\makeatletter
\def\ps@pprintTitle{%
 \let\@oddhead\@empty
 \let\@evenhead\@empty
 \def\@oddfoot{\centerline{\thepage}}%
 \let\@evenfoot\@oddfoot}
\makeatother

\newcommand{\usepackageifexists}[1]{%
     \IfFileExists{#1.sty}{\usepackage{#1}}%
        {\GenericInfo{taglia}{Il package #1 non esiste.}}}

\makeatletter
\def\paragraph{\@startsection{paragraph}{4}%
  \z@\z@{-\fontdimen2\font}%
  {\normalfont\bfseries}}
\makeatother

\usepackage{amsmath,amsthm,amssymb}
\usepackage{amsthm}
\theoremstyle{plain}
\newtheorem{thm}{Theorem}[section]

\newtheorem{lem}[thm]{Lemma}
\newtheorem{conj}[thm]{Conjecture}

\newtheorem{prop}[thm]{Proposition}
\theoremstyle{definition}

\theoremstyle{remark}
\newtheorem{remark}[thm]{Remark}

\newtheorem{ex}[thm]{Example}

\usepackage{upref}
\usepackage{graphicx, caption, subcaption}
\usepackage{color}

\usepackage{algpseudocode}

\usepackage{listings}
\usepackage{epstopdf}
\usepackage{rotating}
\usepackage{mathrsfs}
\usepackage{multirow}
\providecommand{\mathscr}[1]{\mathcal{#1}}

\usepackage{pdfpages}
\usepackage{graphicx}
\usepackageifexists{pdfsync}
\graphicspath{ {figures/} }
\usepackage{booktabs}
\usepackage{hyperref}
\usepackage{geometry, natbib} 
\usepackage{adjustbox}
\usepackage[autostyle]{csquotes}

\def\epsilon{\varepsilon}

\def\R{{\Bbb R}}

\def\N{{\Bbb N}}

\def\1{{\bf 1}}

\def\phi{\varphi}

\def\S{{\cal S}}

\newcommand{\card}{{\rm Card}\;}

\usepackage[toc,page]{appendix}

\def\correspondingauthor{\footnote{Corresponding author. E-mail: felice.iavernaro@uniba.it}}
\begin{document}

\title{A probabilistic approach to the twin prime \\and cousin prime conjectures}

\author[D.Bufalo]{Daniele Bufalo} 
\address[D.Bufalo]{ Dipartimento di Informatica, Universit\`a degli Studi di Bari Aldo Moro, Italy} 
 
\author[M.Bufalo]{Michele Bufalo} 
\address[M.Bufalo]{Dipartimento di Metodi e Modelli per l'Economia, il Territorio e la Finanza,  Universit\`a degli Studi di Roma ``La Sapienza'', Italy} 
\author[Iavernaro]{Felice Iavernaro\correspondingauthor{}} 
\address[Iavernaro]{Dipartimento di Matematica, Universit\`a degli Studi di Bari Aldo Moro, Italy}


\begin{abstract}
We address the question of the infinitude of twin and cousin prime pairs from a probabilistic perspective. Our approach,  based on the novel algorithm introduced in \cite{Bufalos}, partitions the set of integer numbers greater than $2$ in finite intervals of the form $[p_{n-1}^2,p_n^2)$,  $p_{n-1}$ and $p_n$ being two consecutive primes, and evaluates the probability $q_n$ that such an interval contains a twin prime and a cousin prime. Combining Merten's third theorem with the properties of the binomial distribution, we show that $q_n$ approaches $1$ as $n \to \infty$. A study of the convergence properties of the sequence $\{q_n\}$  allows us to propose a new, more stringent conjecture concerning the existence of infinitely many twin and cousin primes. In accord with the Hardy-Littlewood conjecture, it is also shown that twin and cousin primes share the same asymptotic distribution.
\end{abstract}


\maketitle
\section{Introduction}
In the world of mathematics, the greatest unsolved problems often live in number theory, the Riemann hypothesis being one of the most striking  (see, e.g. \cite{Bombieri}). Another prominent example concerns the infinitude of twin primes, namely prime pairs $(p,p+2)$. This is  a commonly accepted statement whose proof has resisted any attack from its first formulation and remained unsolved to the present day. Though seldom ascribed to Euclid, who actually showed the existence of infinitely many prime numbers, the twin prime conjecture implicitly appears in 1849  as a special case of de Polignac's conjecture, which states that any even number can be expressed  as the difference between two consecutive primes in infinitely many ways \cite{Polignac}.  In 1878, J.W.L. Glaisher gave a formal and explicit formulation of the twin prime conjecture thus challenging the mathematical community to prove it (see \cite{Glaisher,Dunham}).

The twin prime conjecture, as well as  other related conjectures, have received support through a wide number of numerical simulations (see, e.g. \cite{Guy}). Starting form 2007, two distributed computing projects,  {\em Twin Prime Search} and {\em PrimeGrid}, have computed several record-largest twin primes.  At the time of this writing, the current largest twin prime pair known is $2996863034895 \cdot 2^{1290000} \pm 1$  which contains $388342$ decimal digits.\footnote{See \url{http://www.primegrid.com/download/twin-1290000.pdf} for the official announcment.} However, searching for largest (twin) prime numbers is always more complicated, due to the information technology limitiations (at least until the quantum memories and quantum computers become reality).

A generalization of this conjecture, known as  the  \emph{Hardy-Littlewood conjecture},  concerns  the distribution of twin primes. Given $x\in \R_+$ with $x\geq 3$, let $\pi_2(x)$  be the enumerative function of twin primes, i.e., the number of primes $p\leq x$ such that $p+2$ is also prime. Then the Hardy-Littlewood conjecture implies that
$$
\pi_2(x)\sim 2C_2 \int_2^x \frac{\mathrm{d}t}{\ln^2 t},
$$    
in the sense that the quotient of the two expressions tends to 1 as $x\to +\infty$. The number $C_2=0.6601618158...$ is called twin prime constant or the Hardy-Littlewood constant.\\
Interesting improvements  have been provided by the Norwegian mathematician Viggo Brun. In 1915, in his work \cite{Brun}, he showed that the sum of the reciprocals of the twin primes converges to a value called  Brun's constant (that is equals to 1.9021605831..., approximately). Should the series have been divergent, as is the case with the sum of the reciprocals of the primes, one would deduce the infinitude of twin primes. In contrast, Brum's result suggests that the twin primes exist in finite number or that, if they are infinite, must represent an infinitesimal fraction of primes.

In the spirit of Brun's idea, sieve methods  may  also help to find certain upper bounds for counting  twin primes. In this context, the strongest result about twin prime conjecture  is provided by Chen's theorem  which states that there are infinitely many primes $p$ such that $(p+2)$ has at most $2$ prime factors \cite{Chen}. Unfortunately, sieve methods often suffer from many limitations, so that new ideas and approaches are needed.

Among the several results weaker than the twin prime conjecture, we recall the one of P. Erd\"os  which states that there exists a constant $c<1$ such that two consecutive primes $(p_n,p_{n+1})$ are less than $c\cdot\log p_n$ \cite{Erdos}. In 2004, D. Goldston et al. (see \cite{Gold1}) reduced such a constant $c$ further to the value $0.08576\dots$ while, in 2009 (see \cite{Gold2}), they established that $c$ can be chosen arbitrarily small, i.e.,
$$
\liminf_{n\to +\infty} \frac{p_n-p_{n-1}}{\log p_{n-1}}=0.
$$
Actually, on April 2013 (see \cite{Zhang}), Y. Zhang proved that there are infinitely many pairs of consecutive primes that differ by a number less than $7\cdot 10^7$, i.e.
$$
\liminf_{n\to +\infty}  (p_n-p_{n-1})<7\cdot 10^7.
$$
On April 2014, this result was improved  by J. Maynard (see \cite{Maynard}) who proved that there exist infinitely many consecutive primes with gap equal to 246. 
Further, assuming the Elliott-Halberstam conjecture and its generalized form, the {\em Polymath} project (see, e.g., \cite{Poly1}, \cite{Poly2}) states that the bound has been reduced to 12 and 6, respectively. Evidently, Zhang's approach is a candidate to enable progress towards a formal proof of the infinitude of twin primes, aiming at showing that the bound above could be reduced to 2.

Moreover, among new ideas helpful to solve the twin prime conjecture we recall those about large gaps between primes.  In 2016, K. Ford et al. (see \cite{Tao1}, \cite{Tao2}) answered  the  Erd\"os question and computed a lower bound for the function counting the largest gap between two consecutive primes below a given integer.

Finally, many probabilistic results can be found in the recent literature, concerning the distribution of primes, such as, for example,  \cite{Baker},\cite{Banks},\cite{Castillo}, \cite{Freiberg}, \cite{Maier}, \cite{Vatwani}, \cite{Tao2011}.
For a detailed review of twin primes the reader can see the survey of J. Maynard \cite{Maynard2}.  

In the present work,  by exploiting  the new algorithm introduced in \cite{Bufalos} and described in the next section, we  give a probabilistic evidence of the twin prime conjecture. In more detail, denoting by $p_n$ the $n$th prime number,   we show that the probability $q_n$  of finding a twin prime inside the interval $ I_{n-1}=[p^2_{n-1},p^2_n)$ tends to 1 as $n \to \infty$. This result, presented in Section \ref{Sec3}, is obtained by combining Merten's third theorem and the  binomial distribution probability function. It turns out that the argument and the tools involved to investigate twin primes remain precisely the same when referred to cousin primes, couples of primes of the form $(p,p+4)$. This circumstance allows us to easily show that these two classes of prime numbers share the same asymptotic distribution, consistently with the Hardy-Littlewood conjecture.  An analytical study of the rate of convergence of $q_n$ to $1$ is then carried out in Section \ref{Sec4}. It makes us confident  that any interval $I_{n-1}$ may indeed  contain twin and cousin primes, which leads us to propose  new, more stringent conjectures at the end of the paper. Some concluding remarks are finally drawn in Section \ref{Sec5}.

\section{A new algorithm for prime numbers}\label{Sec2}
\noindent In this  section, we describe the algorithm  recently introduced in \cite{Bufalos}. It computes all primes progressively by exploring the configuration of a sequence of vectors $C^{(i)}$, $i=3,5,7,\dots,2k+1$ with $k\in \N^\ast$, whose length $l$, in turn, does not remain constant but may increase with $i$. The $j$th prime, when detected, will be added as the $j$th element of a list $P$. Initially, for $k=1$, we set $C^{(3)}=(3)$ and $P=(3)$ (the first prime $2$ is deliberately excluded). The sequence $C^{(i)}$ is recursively defined as follows. Let $l$ be the length of $C^{(i-2)}$.\footnote{In order to keep the notation simpler, we avoid to make the dependence of the length $l$ on the iteration index $i$ explicit.} We first decrease each element of $C^{(i-2)}$ by one unit, setting
$$
\widetilde C^{(i)} = C^{(i-2)}-(1,\dots,1).
$$
We observe that $i$ is a prime number if and only if $\widetilde C^{(i)}$ does not contain null components: in such a case we extend the list $P$ by including $i$ as the new latest element.  $C^{(i)}$ is defined by replacing any zero element in $\widetilde C^{(i)}$ with the corresponding prime number in the list $P$, that is,  for $j=1,\dots,l$,
\begin{equation}
\label{operations}
C^{(i)}(j)= \left\{
\begin{array}{ll}
\widetilde C^{(i)}(j), & \mbox{if } \widetilde C^{(i)}(j)\not= 0,\\
P(j),& \mbox{otherwise}.
\end{array}
\right.
\end{equation}
Moreover, if 
\begin{equation}\label{square0}
\widetilde C^{(i)}(l)=0, \quad \mbox{and} \quad \widetilde C^{(i)}(j)\not=0,\,\forall j<l,
\end{equation}
we increase the length of $C^{(i)}$ by appending the new entry $C^{(i)}(l+1)$ defined as
\begin{equation}\label{square}
C^{(i)}(l+1)=P(l+1)-\frac{d^2}{2} = P(l+1)-\frac{1}{2}(i-P(l+1))\mod P(l+1),
\end{equation}
where $d=P(l+1)-P(l)$ is the $l^{th}$-prime gap.

The presence of a null element in the $j$th position of the vector $C^{(i)}$ means that $P(j)$ divides $i$ (so $i$ cannot be a prime number). In particular, condition \eqref{square0} is fulfilled if and only if  $i=P(l)^2$. In fact, \eqref{square0} means that $P(l)$ divides (but is different from) $i$ and all primes less than $P(l)$ are coprimes with $i$. In this event, starting from the subsequent indices $i+2,i+4,\dots$, the vector $C^{(\cdot)}$ is extended with a new element that brings the information related to the prime number $P(l+1)$. In so doing, the algorithm will correctly detect $i=P(l+1)^2$ as a non-prime number. Consequently, in the interval  $I_l=[P(l-1)^2,P(l)^2)$ the sequence $C^{(i)}$ has constant length $l$ and the primality test  will  be successfully carried out by the algorithm as long as $i\le P(l)^2$. Such an interval will play  an important role in the sequel.   Assuming, by definition, that it is made up of only odd integers, its cardinality is $(P(l)^2-P(l-1)^2)/2$.        

Summing up, the scheme of the algorithm that detects prime numbers up to the integer $n\ge 3$ is the following:

\medskip

\begin{algorithmic}[1]
\State initialize $i\gets 3$,~ $j\gets 1$,~ $P(j)\gets 3$,~ $C^{(i)}(j)\gets 3$,~ $l\gets 1$; 
\While {$i< n$,}:
\State $i \gets i+2$,~ $null\gets 0$;
\For {$j=1,\dots,l$,}:
\State $C^{(i)}(j) \gets C^{(i-2)}(j)-1$
\If {$C^{(i)}(j)=0$,}:
\State $C^{(i)}(j)\gets P(j)$,~ $null \gets null+1$,~ $jmax=j$;
\EndIf
\EndFor 
\If {$null=0$,}: 
\State $P \gets [P, i]$; ($i$ is a prime)
\ElsIf {$null=1$ and $jmax=l$,}:
\State $d\gets P(l + 1)-P(l)$,~ $C^{(i)}(l+1)\gets P(l+1)-d^2/2$,~ $l \gets l+1$;
\EndIf
\EndWhile
\end{algorithmic}

\medskip

As the index $i$ progresses, each digit $C^{(i)}(j)$ behaves like a discrete oscillator (a clock) with period $P(j)$. After each period, namely when $C^{(i)}(j)=0$, we can conclude that $P(j)$ divides $i$. Instead of using the values $1,\dots,P(j)$, we could use any other sequence of symbols such as, for example $0,\dots,P(j)-1$. This latter choice, that will be adopted in the sequel, allows us to exploit the modular arithmetic and remove the first assignment in step 7 of the above algorithm, replacing step 5  with  
\begin{equation}
\label{mod-ar}
C^{(i)}(j) \gets (C^{(i-2)}(j)-1)\mod P(j).
\end{equation}
Table \ref{T1} shows the  sequences $C^{(i)}$ and the vector of primes $P$ resulting as the output of the described procedure with the variant of formula \eqref{mod-ar}. 
\begin{table} \tiny
\begin{center}
\begin{tabular}{c|cccccccc}
\hline\\[-.15cm]
$i:$       & $3$      & $5$      & $7$      &  $9$     &  $11$      & $13$     & $15$      & $17$ \\[.15cm] 
$C^{(i)}:$ & $(3)$    & $(2)$    & $(1)$    &  $(0,3)$ &  $(2,2)$   & $(1,1)$  & $(0,0)$   & $(2,4)$ \\[.15cm]
$P:$       & $P(1)=3$ & $P(2)=5$ & $P(3)=7$ &  $-$     &  $P(4)=11$ & $P(5)=13$& $-$       & $P(6)=17$ \\[.15cm]
\hline\\[-.15cm]
$i:$       & $19$      & $21$    & $23$      & $25$    & $27$    & $29$     & $31$       & $33$    \\[.15cm] 
$C^{(i)}:$ & $(1,3)$   & $(0,2)$ & $(2,1)$   & $(1,0,5)$ & $(0,4,4)$ & $(2,3,4)$  & $(1,2,4)$    & $(0,1,1)$ \\[.15cm]
$P:$       & $P(7)=19$ & $-$     & $P(8)=23$ & $-$     & $-$     & $P(9)=29$& $P(10)=31$ & $-$     \\[.15cm]
\hline\\[-.15cm]
$i:$       & $35$      & $37$        & $39$      & $41$    & $43$    & $45$     & $47$       & $49$    \\[.15cm] 
$C^{(i)}:$ & $(2,0,0)$   & $(1,4,6)$ & $(0,3,5)$   & $(2,2,4)$ & $(1,1,3)$ & $(0,0,2)$  & $(2,4,1)$    & $(1,3,0,3)$ \\[.15cm]
$P:$       & $-$ & $P(11)=37$         & $-$ & $P(12)=41$     & $P(13)=43$     & $-$& $P(14)=47$ & $-$     \\[.15cm]
\hline\\[-.15cm]
$i:$       & $51$      & $53$        & $55$      & $57$    & $59$    & $61$     & $63$       & $65$    \\[.15cm] 
$C^{(i)}:$ & $(0,2,6,2)$   & $(2,1,5,1)$ & $(1,0,4,0)$   & $(0,4,3,10)$ & $(2,3,2,9)$ & $(1,2,1,8)$  & $(0,1,0,7)$    & $(2,0,6,6)$ \\[.15cm]
$P:$       & $-$ & $P(15)=53$         & $-$ & $-$     & $P(16)=59$     & $P(17)=61$& $-$ & $-$     \\[.15cm]
\hline\\[-.15cm]
$i:$       & $67$        & $69$        & $71$      & $73$    & $75$    & $77$     & $79$       & $81$    \\[.15cm] 
$C^{(i)}:$ & $(1,4,5,5)$ & $(0,3,4,4)$ & $(2,2,3,3)$   & $(1,1,2,2)$ & $(0,0,1,1)$ & $(2,4,0,0)$  & $(1,3,6,10)$    & $(0,2,5,9)$ \\[.15cm]
$P:$       & $P(18)=67$ & $-$  & $P(19)=71$ & $P(20)=73$ & $-$ & $-$& $P(21)=79$ & $-$     \\[.15cm]
\hline\\[-.15cm]
$i:$       & $83$        & $85$        & $87$      & $89$    & $91$    & $93$     & $95$       & $97$    \\[.15cm] 
$C^{(i)}:$ & $(2,1,4,8)$ & $(1,0,3,7)$ & $(0,4,2,6)$   & $(2,3,1,5)$ & $(1,2,0,4)$ & $(0,1,6,3)$  & $(2,0,5,2)$    & $(1,4,4,1)$ \\[.15cm]
$P:$       & $P(22)=83$ & $-$  & $-$ & $P(23)=89$ & $-$ & $-$& $-$ & $P(24)=97$     \\[.15cm]
\hline\\[-.15cm]
$i:$       & $99$        & $101$        & $103$      & $105$    & $107$    & $109$     & $111$       & $113$    \\[.15cm] 
$C^{(i)}:$ & $(0,3,3,0)$ & $(2,2,2,10)$ & $(1,1,1,9)$   & $(0,0,0,8)$ & $(2,4,6,7)$ & $(1,3,5,6)$  & $(0,2,4,5)$    & $(2,1,3,4)$ \\[.15cm]
$P:$       & $-$ & $P(25)=101$  & $P(26)=103$ & $-$ & $P(27)=107$ & $P(28)=109$ & $-$ & $P(29)=113$     \\[.15cm]
\hline\\[-.15cm]
$i:$       & $115$        & $117$        & $119$      & $121$        \\[.15cm] 
$C^{(i)}:$ & $(1,0,2,3)$ & $(0,4,1,2)$ & $(2,3,0,1)$   & $(1,2,6,0,11)$ \\[.15cm]
$P:$       & $-$ & $-$  & $-$ & $-$      \\[.15cm]
\hline
\end{tabular}
\end{center}
\caption{Algorithm steps up to $100$. The scheme shows the computation of  the sequence $C^{(i)}$ and the vector $P$ for all odd indices $i\le 121 = 11^2$.}
\label{T1}
\end{table}
\begin{remark}\label{R1} \em
Let $l$ be the size of $C^{(\cdot)}$ that, as was already observed,  remains constant in the set $I_l:=[P^2(l-1),P^2(l))$. We have than that every number $x\in I_l$ satisfies the following system
\begin{equation}
\begin{cases}
x+2C^{(x)}(1)\equiv 0\qquad mod\,P(1) \\ x+2C^{(x)}(2)\equiv 0\qquad mod\,P(2)\\ \vdots \\ x+2C^{(x)}(l)\equiv 0\qquad mod\,P(l);
\end{cases}
\end{equation}
whose solution, in turn, have to satisfy the constrain $x\in I_l$. Hence, there exists a biunivocal relation between $x\in I_l$ and $C^{(x)}$ (between $I_l$ and the sequences $C^{(x)}$ of length $l$).
\end{remark}

\section{On the infinitude of twin and cousin primes}
\label{Sec3}
\noindent From the algorithm explained in Section \ref{Sec2}, the following characterization of primes,  twin primes and cousin primes in terms of the sequence $C^{(i)}$ is readily deduced.
\begin{prop}[Primes, twin primes and cousin primes characterization]
\label{prop21}
Let  $h\ge 3$ be a given integer. The following propositions hold true:
\begin{itemize}
\item[(a)]  $h$ is a prime number if and only if the tuple $C^{(h)}$ does not contain zeroes, i.e.,
\begin{equation}
\label{relazione_primi}
C^{(h)}(j)\not= 0, \qquad \forall j=1,\dots, l,
\end{equation}
$l$ being the length of $C^{(h)}$.
\item[(b)] $(h,h+2)$ is a pair of twin primes if and only if the tuple $C^{(h)}$ does not contain neither $0$ nor $1$, i.e.,
\begin{equation}
\label{relazione_gemelli}
C^{(h)}(j)\not= 0,1, \qquad \forall j=1,\dots, l.
\end{equation}
\item[(c)] $(h,h+4)$ is a pair of cousin primes if and only if the tuple $C^{(h)}$ does not contain neither $0$ nor $2$, i.e.,
\begin{equation}
\label{relazione_cugini}
C^{(h)}(j)\not= 0,2, \qquad \forall j=1,\dots, l.
\end{equation}
\end{itemize}
\end{prop} 

We notice that (\ref{relazione_gemelli})  and (\ref{relazione_cugini}) are similar in nature: both require the absence of two distinct digits in the tuple $C^{(h)}$. The argument introduced below makes sense regardless of the choice of these two digits and, consequently, the results we will derive hold true for both twin and cousin prime pairs. For this reason, for the sake of simplicity, we address the discussion to the former subset of primes and then export the results to the latter.

In the following, $\mathcal{P}=\{ 2,3,5,7,11,...\}$ defines the sequence of all primes and its elements will be denoted $p_i$, $i=1,\dots$. Clearly, apart from the very first element $2$, $\mathcal{P}$ matches the vector $P$ resulting from running infinitely many iterations of the algorithm in Section \ref{Sec2}, and hence $P(i)=p_{i+1}$, for all indices $i\ge 1$.\footnote{Dealing with the set $\mathcal{P}$ is more convenient to conform with the standard notation in the literature.} Furthermore, as is usual,  the notation $A(x) \sim B(x)$, for real or integer $x$, means that $\lim_{x\rightarrow +\infty}A(x)/B(x)=1$, so both quantities have the same asymptotic behavior. 
\begin{lem}\label{lemma1}
For all integrals $n\ge 2$, the following relation holds true:
\begin{equation}\label{rel1}
\prod_{i= 2}^{n} \biggl( 1-\frac{2}{p_i}\biggr)=C_2(n)\biggl( \prod_{i= 2}^{n} \biggl( 1-\frac{1}{p_i}\biggr)\biggr)^2=4C_2(n)\biggl( \prod_{i= 1}^{n} \biggl( 1-\frac{1}{p_i}\biggr)\biggr)^2,
\end{equation} 
where 
\begin{equation}\label{C2n}
C_2(n)=\prod_{i=2}^{n} \frac{p_i(p_i-2)}{(p_i-1)^2}=\prod_{i=2}^{n} \biggl( 1-\frac{1}{(p_i-1)^2}\biggr).
\end{equation}
Consequently, the following asymptotic estimation holds true:
\begin{equation}\label{rel1a}
\prod_{i= 2}^{n} \biggl( 1-\frac{2}{p_i}\biggr)\sim C_2\biggl( \prod_{i= 2}^{n} \biggl( 1-\frac{1}{p_i}\biggr)\biggr)^2=4C_2\biggl( \prod_{i= 1}^{n} \biggl( 1-\frac{1}{p_i}\biggr)\biggr)^2,
\end{equation} 
where $C_2$ is the Hardy-Littlewood constant
\begin{equation}\label{C2}
C_2=\lim_{n\rightarrow +\infty}C_2(n)=\prod_{i=2}^{+\infty} \frac{p_i(p_i-2)}{(p_i-1)^2}=\prod_{i=2}^{+\infty} \biggl( 1-\frac{1}{(p_i-1)^2}\biggr)=0.6601618\cdots .
\end{equation}
\end{lem} 
\begin{proof}
For any integer $x \ge 2$, a direct computation shows that
\begin{equation}
\label{simple-rel}
1-\frac{2}{x}= \biggl( 1-\frac{1}{x}\biggr)^2\biggl( 1-\frac{1}{(x-1)^2}\biggr).
\end{equation}
To prove (\ref{rel1}), we exploit an induction argument on $n$. When $n=2$, we see that (\ref{rel1}) is nothing but (\ref{simple-rel}) with $x=p_2$.
Assuming that the statement holds true for $n$, we prove it  for $n+1$. Exploiting again (\ref{simple-rel}), we get
$$ 
\begin{array}{rcl} \displaystyle
\prod_{i= 2}^{n+1} \left( 1-\frac{2}{p_i}\right)&=& \displaystyle \biggl( 1-\frac{2}{p_{n+1}}\biggr) \prod_{i= 2}^{n} \biggl( 1-\frac{2}{p_i}\biggr) \\[.5cm]
&=&\displaystyle \biggl( 1-\frac{1}{p_{n+1}}\biggr)^2\biggl( 1-\frac{1}{(p_{n+1}-1)^2}\biggr) \prod_{i= 2}^{n} \biggl( 1-\frac{1}{p_i}\biggr)^2\biggl( 1-\frac{1}{(p_i-1)^2}\biggr) \\[.5cm]
&=&\displaystyle \prod_{i= 2}^{n+1} \biggl( 1-\frac{1}{p_i}\biggr)^2\biggl( 1-\frac{1}{(p_i-1)^2}\biggr).
\end{array}
$$
Relation (\ref{rel1a}) is a direct consequence of (\ref{rel1}), after noticing that  $\lim_{n\rightarrow +\infty} C_2(n)=C_2$ (see (\ref{C2})).
\end{proof} 
Combining Lemma \ref{lemma1} with Mertens' third theorem \cite{Tenenbaum}
$$
\lim_{n\rightarrow +\infty} \log p_n \prod_{i= 1}^{n} \biggl( 1-\frac{1}{p_i}\biggr) = e^{-\gamma},
$$
where $\gamma \simeq 0.57721$ is the  Euler-Mascheroni constant, the following result immediately follows.
\begin{lem}\label{lemma2}
The following asymptotic estimation holds true:
\begin{equation}\label{rel1b}
\prod_{i= 2}^{n} \biggl( 1-\frac{2}{p_i}\biggr)\sim 4C_2  e^{-2\gamma} \frac{1}{\log^2 p_n}.
\end{equation} 

\end{lem}

Now, for $n\ge 3$, let us consider the sequence $C^{(k)}$, with  $k\in I_{n-1}=[P(n-2)^2,P(n-1)^2)=[p_{n-1}^2,p_{n}^2)$. As was emphasized in Section \ref{Sec2}, we recall that
\begin{itemize}
\item[(i)] all these  $C^{(k)}$ share the same length $n-1$ since they only involve the primes $p_2,\dots,p_{n}$;
\item[(ii)] there are $T_n:=\card(I_{n-1})=(p_n^2-p_{n-1}^2)/2$ such tuples;
\item[(iii)] each element $C^{(k)}(j-1)$, $j=2,\dots,n$, has periodicity $p_{j}$ (a prime number) with respect to the index $k$.
\end{itemize}
Due to property $(iii)$, each element $C^{(k)}(j-1)$, $j=2,\dots,n$, may be thought of as an event in the sample space $\Omega_j$ composed of all the possible $p_j$ outcomes $0,1,\dots,p_j-1$. If all the events in $\Omega_j$ were equiprobable, we could conclude that, picking a random tuple $C^{(k)}$ in $I_{n-1}$,  the probability that  $C^{(k)}(j-1)\not = 0,1$ would be equal to $(1-2/p_j)$.\footnote{Of course, the same  holds true for the probability that $C^{(k)}(j-1)\not = 0,2$.}   Assuming, for a while, that two events in $\Omega_{j_1}$ and $\Omega_{j_2}$ are independent,  we could conclude (see Proposition \ref{prop21}) that the probability that  such a tuple $C^{(k)}$ identifies a twin pair is just $\prod_{i= 2}^{n} ( 1-2/p_i)$, which would justify the estimations in Lemmas  \ref{lemma1} and \ref{lemma2}. 
Before introducing the main result of our study, a discussion about the extent to which the equiprobability and independency assumptions are retained by the subset 
\begin{equation}
\label{Sn}
\S_{n-1}:=\{C^{(k)}: ~k\in I_{n-1}\} \subset \Omega_2\times \cdots \times \Omega_n
\end{equation}
is presented in the following two subsections.
\subsection{Equiprobability}
For the counter $C^{(k)}(j-1)$, the assumption of equiprobability is valid only in the event that $T_n$ is a multiple of $p_j$ which, in general, is not the case. On the other hand, if $p_j$ is much smaller than $T_n$, the actual probability that $C^{(k)}(j)\not = 0,1$ gets very close to $(1-2/p_j)$. Before formalizing this concept we consider a preliminary example.
\begin{ex} 
\label{example1}
\em 
Consider  the $36$ tuples $C^{(k)}$ of length $4$ in the interval $I_4 =[7^2,11^2)$ displayed in Table \ref{T1}. The first two columns of Table \ref{T2} compare the relative frequencies of the outcomes $C^{(k)}(j)\not = 0,1$ (actual probabilities) with the reference probabilities  $(1-2/p_j)$ one would obtain when $p_j$ divides $T_n$. For increasing values of $n$ and thus of $T_n$, the observed discrepancies of the behavior of these first $4$ digits become smaller and smaller. As an example, Table \ref{T2} displays the results for  $I_{14}=[43^2,47^2)$ and $I_{29}=[109^2,113^2)$. Notice that there is no error in $C^{(k)}(1)$ since $T_n$ is always a multiple of $p_2=3$. In fact, from Fermat's little theorem, namely $a^{p-1}-1$ is a multiple of $p$ for any pair $(a,p)$ of coprime integers, setting $p=3$ and $a=p_n$, we see that 
$$
T_n=\frac{p_n^2-p_{n-1}^2}{2}= \frac{(p_n^2-1)-(p_{n-1}^2-1)}{2}  
$$
is divisible by $3$.
\begin{table} 
\begin{center}
\begin{tabular}{c|cccc}
            & $(1-2/p_j)$ & \footnotesize \begin{tabular}{c}  relative frequency \\ in $I_4=[7^2,11^2)$ \end{tabular} & \footnotesize \begin{tabular}{c}  relative frequency \\ in $I_{14}=[43^2,47^2)$ \end{tabular} & \footnotesize \begin{tabular}{c}  relative frequency \\ in $I_{29}=[109^2,113^2)$ \end{tabular}\\ 
\hline\\[-.15cm]
$C{(k)}(1)$ &  $1/3$      & $1/3$       &  $1/3$      & $1/3$    \\[.15cm] 
$C{(k)}(2)$ &  $0.6$      & $0.611111$  &  $0.602564$ & $0.599099$    \\[.15cm]
$C{(k)}(3)$ &  $0.714285$ & $0.694444$  &  $0.717948$ & $0.716216$    \\[.15cm] 
$C{(k)}(4)$ &  $0.818181$ & $0.805555$  &  $0.820512$ & $0.819819$    \\[.15cm]  
\hline
\end{tabular}
\end{center}
\caption{Comparison between the  reference and actual probabilities of the event $C^{(k)}(j)\not = 0,1$ in the sets $I_4$, $I_{12}$ and $I_{29}$. The numbers are rounded after the sixth digit.}
\label{T2}
\end{table}
\end{ex}
To account for this anomaly, for $j=2,\dots,n$, we introduce the relative frequencies
\begin{equation}
\label{qj}
q(j,T_n) = \frac{\card(\{C^{(k)}(j-1)\not = 0,1\}_{k\in I_{n-1}})}{T_n},
\end{equation}
and wish to relate them with the corresponding reference values $1-2/p_j$ which would occur under the equiprobability condition.\footnote{Actually $q(2,T_n)=1-1/p_2$, as was explained at the end of Example \ref{example1}.} 
\begin{lem}\label{lemmaTheta} The following estimation holds true for $j=2,\dots,n$:
\begin{equation}
\label{qjest}
q(j,T_n) = \left( 1-\frac{2}{p_j}\right)\left( 1+ O\Big(\frac{1}{n\log n}\Big)\right).  
\end{equation}
Consequently, 
\begin{equation}
\label{qjprod}
\prod_{j=2}^n q(j,T_n) = \vartheta_n \prod_{j=2}^n  \left(1-\frac{2}{p_j}\right),  \mbox{with } \lim_{n \rightarrow +\infty} \vartheta_n = \lim_{n \rightarrow +\infty} \left( 1+ O\Big(\frac{1}{n\log n}\Big) \right)^{n-1} = 1.
\end{equation}
\end{lem}
\begin{proof}
For $j=2,\dots,n$, let
$$
T_n=\nu_j p_j +r_j, \quad \mbox{with} ~ r_j<p_j.
$$
This means that the $(j-1)$th counter performs $\nu_j$ complete periods plus $r_j$ further steps. 
According to whether or not $0,1$ occur in these extra steps, we bound  $q(j,T_n)$ as
\begin{equation}
\label{luq}
1-2\frac{\nu_j+1}{T_n} \le q(j,T_n) \le 1-2\frac{\nu_j}{T_n}. 
\end{equation}
For the upper bound, we have
$$
\begin{array}{rcl}
\displaystyle  1-2\frac{\nu_j}{T_n} &=& \displaystyle 1-\frac{2}{p_j+\frac{r_j}{\nu_j}} ~=~ 1-\frac{2}{p_j(1 +\frac{r_j}{\nu_jp_j})} ~=~ 1-\frac{2}{p_j} \left(1 +O\Big(\frac{r_j}{\nu_jp_j}\Big)\right) \\[.5cm]
&=& \displaystyle 1-\frac{2}{p_j} +O\Big(\frac{2 r_j}{\nu_jp_j^2}\Big) ~=~ \left(1-\frac{2}{p_j}\right) \left(1 +O\Big(\frac{2 r_j}{\nu_jp_j^2}\Big)\right).
\end{array}
$$
Furthermore, from $-r_j>-p_j\ge -p_n$ and
$$
T_n=\frac{p_n^2-p_{n-1}^2}{2}\ge \frac{p_n^2-(p_{n}-2)^2}{2}=2p_n-2,
$$ 
the equality being possible when $(p_{n-1},p_n)$ is a prime twin pair, it follows that $T_n-r_j>p_n-2$, and hence we get
$$
\begin{array}{rcl}
\displaystyle  
\frac{2 r_j}{\nu_jp_j^2} &<&  \displaystyle  \frac{2 p_j}{\nu_jp_j^2} ~=~ \frac{2}{\nu_jp_j} ~=~ \frac{2}{T_n -r_j} ~<~ \frac{2}{p_n-2} ~=~ O\Big(\frac{1}{p_n}\Big) ~=~ O\Big(\frac{1}{n \log n}\Big). 
\end{array}
$$
The same result holds true for the lower bound in (\ref{luq}), so that \eqref{qjest} and \eqref{qjprod} are definitely deduced.
\end{proof}

\subsection{Independency}
Concerning the independency assumption, both the specific features of the algorithm and the primality of the period of each digit in  $C^{(k)}$ play a key role. Starting from any arbitrary (n-1)-tuple configuration $\gamma^{(1)}=(d_2,\dots,d_{n})$, with $0\le d_i < p_{i}$,  we consider the sequence
\begin{equation}
\label{algmod}
\gamma^{(i)}= \big(\gamma^{(i-1)}-(1,\dots,1)\big ) \mod (p_2,\dots,p_n), \quad i=2,\dots \Gamma_n=\prod_{j=2}^{n}p_j,
\end{equation}
which is therefore defined by the same algorithm described in Section \ref{Sec2} but iterated for $\Gamma_n$ rather than $T_n$ steps. Since  $p_j$ are prime numbers, the tuples  $\gamma^{(k)}$ are all different from each other and covers the whole sample space $\Omega_2\times \cdots \times \Omega_n$. Evidently, in this larger space, the event
$$
E = ``\gamma^{(k)}(i) \not =0,1, \quad \mbox{for} ~i=1,\dots,n-1"
$$
has an associated probability precisely equal to $\prod_{j=2}^{n}(1-2/p_j)$. 

We notice that, independently of the starting value $\gamma^{(1)}$, by decreasing each digit of $\gamma^{(i-1)}$ of one unit per step (in modular arithmetic), the iteration function defined at (\ref{algmod}) changes all the entries of  the input tuple $\gamma^{(i-1)}$ and, as the steps proceed, it put them in periodic circulation. As a consequence, the $\Gamma_n \prod_{j=2}^{n}(1-2/p_j)$ tuples which do not contain  the integers $0$ and $1$ are spread along the sequence at the same manner as are the tuples that, for example, do not contain $0$ and $2$ or $1$ and $2$ respectively. In other words, in terms of distributions of the digits in the tuples $\gamma^{(k)}$ there is no bias introduced by the algorithm, so that any two random subsequences of the same length  share the same probability that a tuple inside them is free of $0$ and $1$, thus detecting a twin prime (or $0$ and $2$ for cousin primes).
 
    Now, the sequence $C^{(k)}$, $k=p_{n-1}^2,\dots,p_n^2$, is just a section of length $T_n$ of the $\gamma^{(i)}$'s and, since $T_n<\Gamma_n$, we may think of the set $S_{n-1}$ containing the tuples  $C^{(k)}$ (see \eqref{Sn}) as a random {\em sample set} extracted from the {\em population set }  $\Omega_2\times \cdots\times \Omega_{n}$. We therefore expect that the relative distribution of digits in the sample  remains the same as that in the population,  modulo a noise due to the non completely uniform distribution along the sequence $\gamma^{(i)}$ of  the tuples belonging to $E$, the equiprobability argument discussed in Lemma \ref{lemmaTheta} and the choice of the initial configuration. To formally state this almost independency property, we may assume that $0<a<b<+\infty$ exist such that the probability $Q_n$ of the event $E$ conditioned to $S_{n-1}$ is 
\begin{equation}
\label{indip-prob}
Q_n = \Theta_n \prod_{j=2}^{n}(1-2/p_j), \quad \mbox{with } a< \Theta_n <b.
\end{equation}    

Below we introduce an illustrative example while, in the next section, we derive an asymptotic estimation of $\Theta_n$.
\begin{ex}
\label{example2}
\em With reference to Table \ref{T1}, we compute the relative distribution of the digits ranging from $0$ to $10$ in all the $36$ tuples of length $4$ (i.e. in $I_4$) and compare it with the corresponding distribution in the complete set $\Omega_2\times \Omega_3 \times \Omega_4 \times \Omega_5$ made up of all the $11\cdot 7\cdot 5 \cdot 3=11500$ configurations of a tuple of length $4$  whose $(j-1)$th entry may range from 0 to $p_{j}-1$, $j=2,\dots,4$. The left picture in Figure \ref{Figdistribution} shows a grey bar plot for the former distribution and a light grey bar plot for the latter. We see that, as a result of the definition of the algorithm, they are very close to each other. As a comparison, in the right plot of   Figure \ref{Figdistribution} we also report the distribution of the digits inside the tuples of length $8$. The Kolmogorov-Smirnov test may be employed to estimate the distance between the sample and reference distribution functions, the null hypothesis being that the sample is drawn from the reference distribution. The $p$-values related to the distributions displayed in Figure \ref{Figdistribution} are $0.37$ and $0.84$ respectively. Furthermore, for increasing values of $n$ the associated $p$-values tend to get very close to $1$; for example, the $p$-values related to the distributions arising from the tuples in $I_{14}$ and $I_{29}$ (see the last  two columns in Table \ref{T2}) are $0.81$ and $0.98$ respectively. Repeating the experiment by allowing different starting tuples in the algorithm, we have experienced similar results. Consequently, the sample set reproduces quite faithfully the reference distribution function, so that we may consider each $C^{(k)}$ to be randomly picked from the population set and apply the independency argument discussed above.       
\begin{figure}[htbp]
\begin{center}
\includegraphics[width=0.49\textwidth]{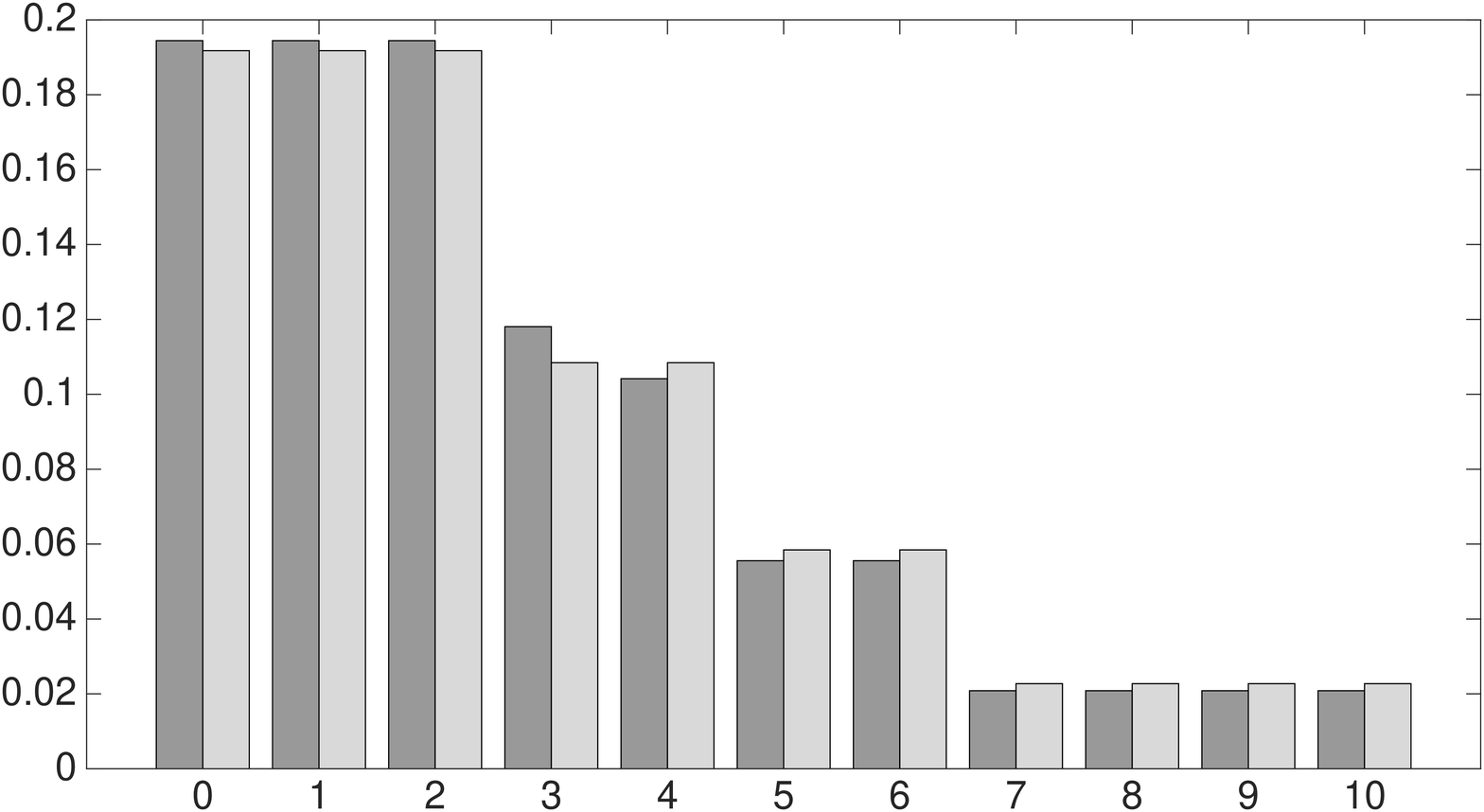}
\includegraphics[width=0.49\textwidth]{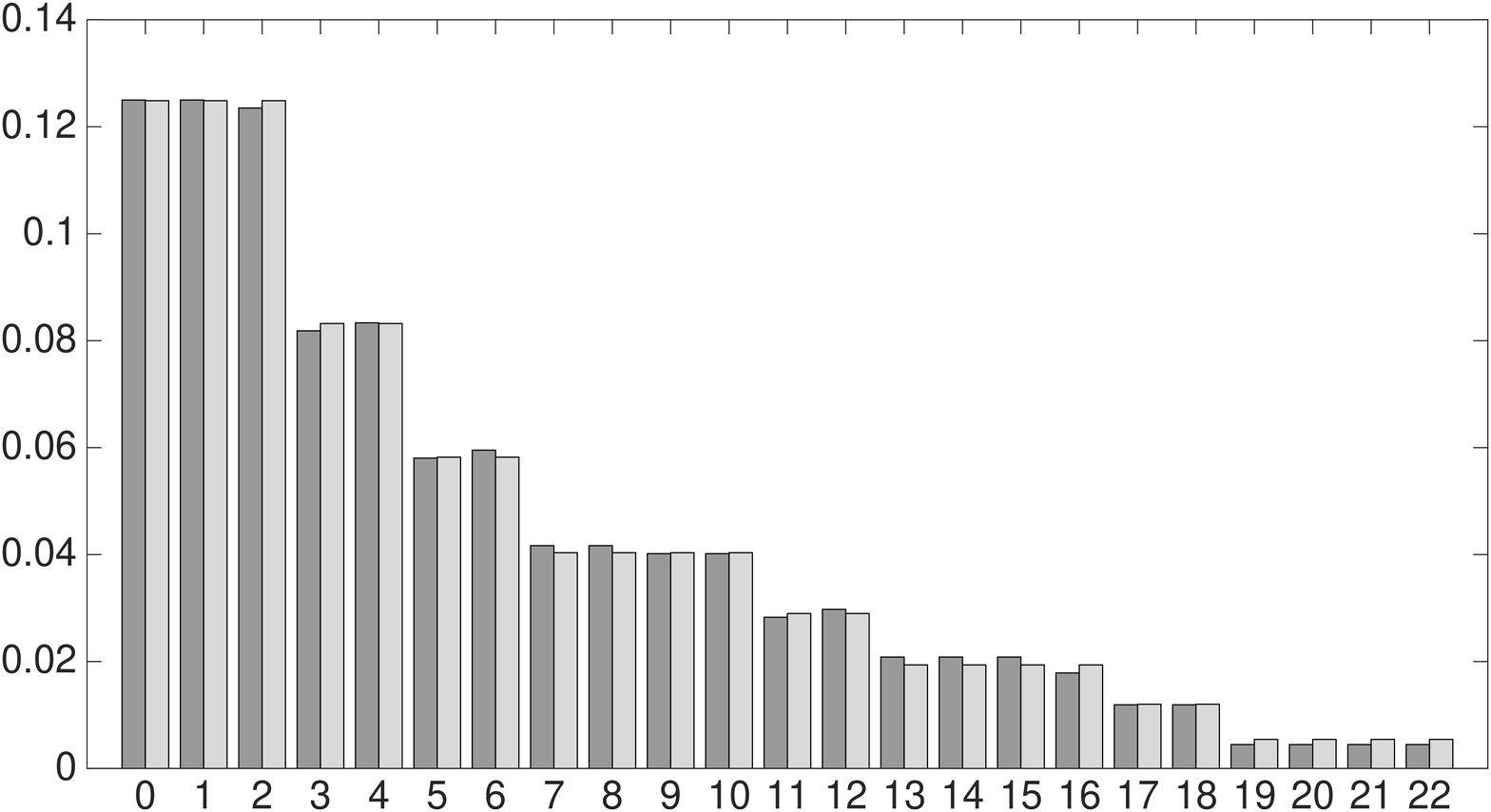}
\end{center}
\caption{Relative frequencies of the digits appearing in the tuples of lengths $4$ (left picture) and $8$ (right picture). The dark grey bars are related to the set of tuples $C^{(k)}$ generated by the algorithm while the light grey bars refer to the set of all the possible tuples of length $4$ and $8$ that can be formed by allowing that their $(j-1)$th element ranges from $0$ to $p_{j}$.}
\label{Figdistribution}
\end{figure}
\end{ex}

\subsection{Main result}
We are now in the right position to state the main result of our study.
\begin{thm}
\label{mainthm}
Under the assumption \eqref{indip-prob}, the probability  that a twin pair lies in the interval $I_{n-1}=[p_{n-1}^2,p_n^2)$ approaches $1$ as $n$ tends to infinity. The same holds true for cousin pairs.   
\end{thm}
\begin{proof}
The probability that no twin/cousin pairs lie in $I_{n-1}$ is $(1-Q_n)^{T_n}$. 
From Lemma \ref{lemma2} we then have
$$ \begin{array}{rcl}
0 &\le& \displaystyle \lim_{n\rightarrow +\infty} (1-Q_n)^{T_n} ~\le~ \lim_{n\rightarrow +\infty} \left(1-\Theta_n \prod_{i= 2}^{n} \left( 1-\frac{2}{p_i}\right)\right)^{2(p_n-1)}\\[.5cm]
& =&\displaystyle  \lim_{n\rightarrow +\infty} \left(1-  4\Theta_n C_2  e^{-2\gamma} \frac{1}{\log^2 p_n}\right)^{2(p_n-1)} ~=~0.
\end{array}
$$
\end{proof}

\section{Novel twin and cousin prime conjectures}
\label{Sec4}
In this section, we wish to get  more insight into the convergence properties of the sequence defined in Theorem  \ref{mainthm}. These features, together with some numerical evidence, will enable us to propose a stronger version of the standard twin and cousin prime conjectures. All numerical tests have been performed in Matlab and Magma environments.
\subsection{Estimation of $\Theta_n$}
\begin{figure}[htbp]
\begin{center}
\hspace*{-.6cm}
\includegraphics[width=0.49\textwidth]{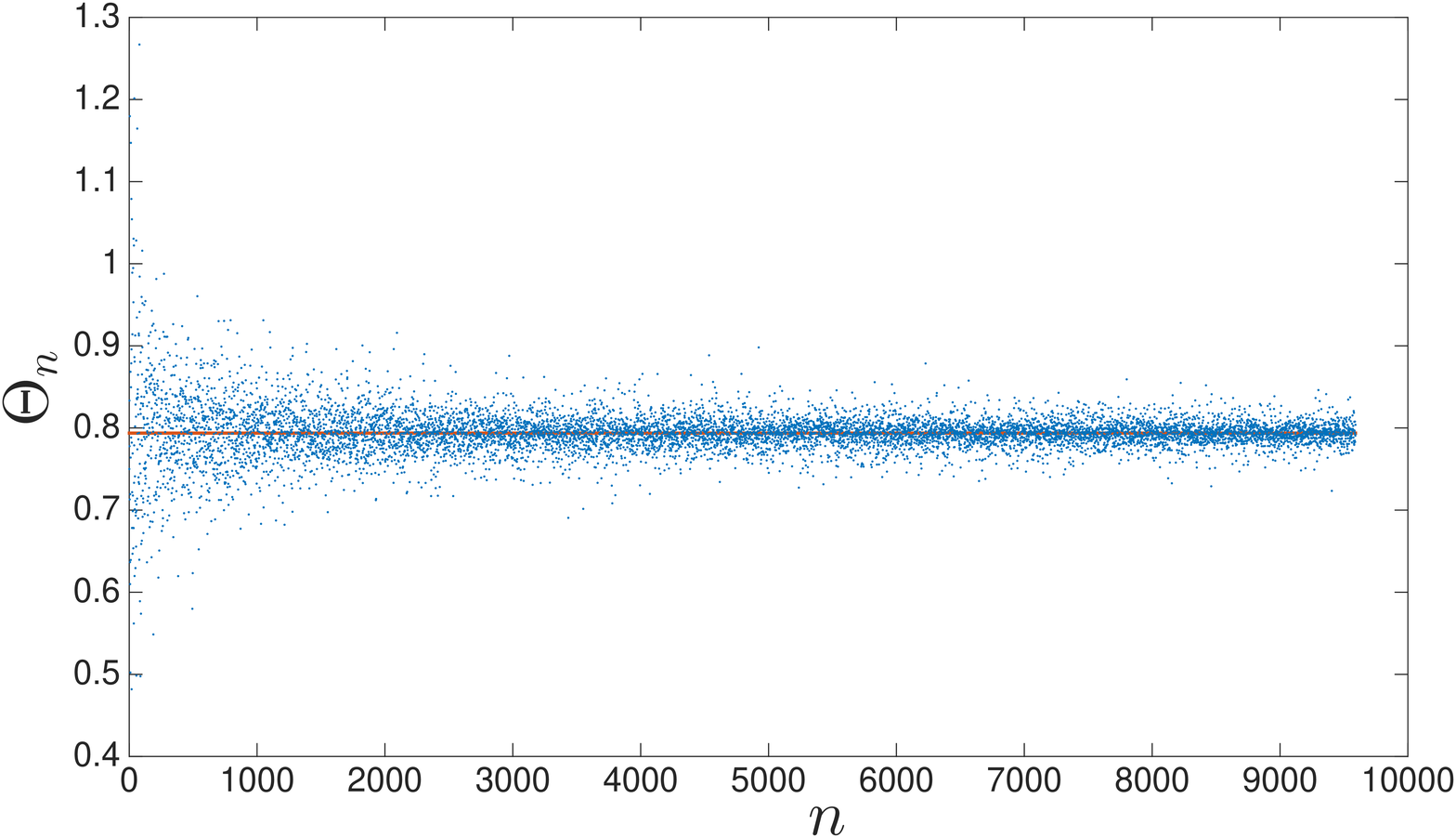}
\includegraphics[width=0.49\textwidth]{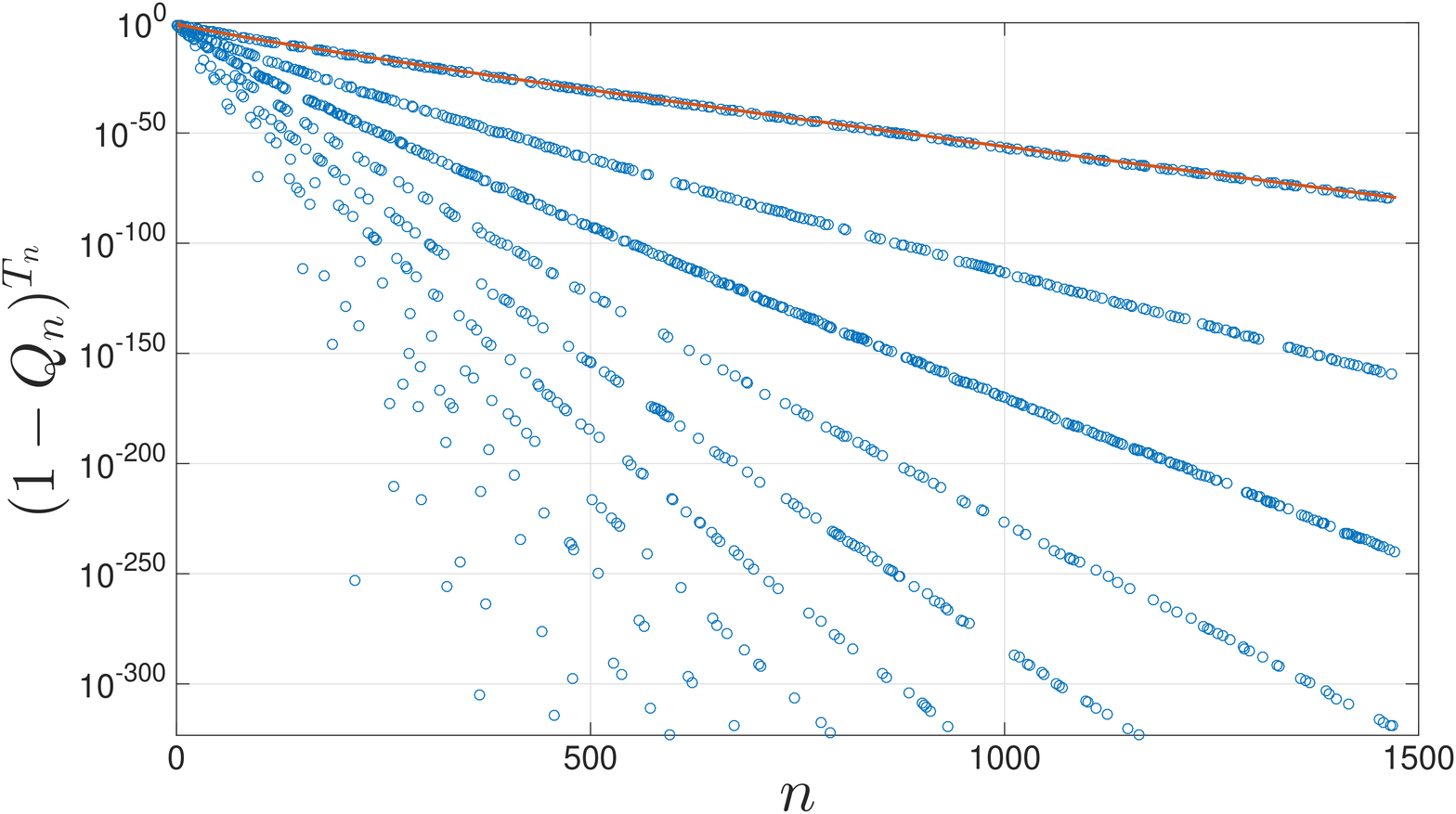}
\end{center}
\caption{Left picture: Values of $\Theta_n$ (dots) and their mean value (solid line). Right picture: decay of the the probability  $(1-Q_n)^{T_n}$ that no twin pair lies in $I_{n-1}$ (small circles) and the analytic approximation of the uppermost branch (solid line).}
\label{FigQn}
\end{figure}

An asymptotic approximation of the sequence $\{\Theta_n\}$ introduced in (\ref{indip-prob}) may be retrieved by making the asymptotic behavior of the probability $Q_n$, namely
\begin{equation}
\label{Qnt}
Q_n \sim 4\Theta_n C_2  e^{-2\gamma} \frac{1}{\log^2 p_n}
\end{equation} 
consistent with the Hardy-Littlewood conjecture
\begin{equation}
\label{HL}
\pi_2(x)\sim 2C_2 Li_2(x),
\end{equation} 
where $\pi_2(x)$ denotes the enumerative function of twin (and cousin) primes up to $x$ and 
\begin{equation}
\label{li2}
Li_2(x)=\int_2^x \frac{\mathrm{d}t}{\log^2 t}.
\end{equation}   
On the basis of formula \eqref{Qnt},  the expected number of twin primes in the interval $I_{n-1}=[p_{n-1},p_n)$, for large values of $n$, is obtained by multiplying $Q_n$ by $T_n$, the cardinality of $I_{n-1}$: 
\begin{equation}
\label{B2I}
T_n Q_n \sim 4 T_n \Theta_n C_2  e^{-2\gamma} \frac{1}{\log^2 p_n}
\end{equation} 
On the other hand, for the same interval $I_{n-1}$, the Hardy-Littlewood estimation would give an expected number of twin primes equal to
\begin{equation}
\label{HL1}
\pi_2(p_n^2) - \pi_2(p_{n-1}^2) \sim 2C_2 \int_{p_{n-1}^2}^{p_n^2} \frac{\mathrm{d}t}{\log^2 t}.
\end{equation}
Since  $\log^2 t$ is strictly increasing for any $t\geq 2$, we have
\begin{equation}
\label{doublebound}
\frac{2T_n}{\log^2 p_n^2}< \int_{p_{n-1}^2}^{p_n^2} \frac{\mathrm{d}t}{\log^2 t} <\frac{2T_n}{\log^2 p_{n-1}^2}.
\end{equation}
A comparison with \eqref{B2I} leads us to investigate the extent to which the leftmost term in \eqref{doublebound} approximates the middle integral term. To this end, also considering that these quantities tend to increase with $n$, we study the relative error
\begin{equation}
\label{En}
E_n = \frac{\int_{p_{n-1}^2}^{p_n^2} \frac{\mathrm{d}t}{\log^2 t} - \frac{2T_n}{\log^2 p_n^2}}{\int_{p_{n-1}^2}^{p_n^2}\frac{\mathrm{d}t}{\log^2 t} }
\end{equation}  
which, besides the order of magnitude, gives us information about the number of significant digits in the approximation. Exploiting \eqref{doublebound} and Bertrand's postulate, i.e. $p_{n}<2p_{n-1}$ for any $n\ge 1$, we obtain
$$
\begin{array}{rcl}
E_n & < & \displaystyle \frac{\frac{2T_n}{\log^2 p_{n-1}^2} - \frac{2T_n}{\log^2 p_n^2}}{\frac{2T_n}{\log^2 p_n^2} } ~ = ~ \frac{\log^2 p_n^2}{\log^2 p_{n-1}^2} -1 ~ < ~ \frac{\log^2 4 p_{n-1}^2}{\log^2 p_{n-1}^2} -1 \\[.6cm]
&=&  \displaystyle \frac{\log^2 4 + 2 \log 4 \log p_{n-1}^2}{\log^2 p_{n-1}^2}~=~ O \left(\frac{1}{\log p_{n-1}}\right).
\end{array}
$$
Thus $\lim_{n \rightarrow +\infty} E_n=0$, which means that (see \eqref{En})  $\frac{2T_n}{\log^2 p_n^2}  \sim \int_{p_{n-1}^2}^{p_n^2} \frac{\mathrm{d}t}{\log^2 t}$. Therefore, \eqref{B2I} and \eqref{HL1} are equivalent in the event that $\lim_{n \rightarrow +\infty} \Theta_n = e^{2\gamma}/4 \simeq 0.79305$. Inserting this estimation in \eqref{Qnt} we conclude that
\begin{equation}
\label{Qnt_lim}
Q_n \sim  \frac{C_2 }{\log^2 p_n}.
\end{equation} 
The left picture of Figure \ref{FigQn} shows the values of $\Theta_n=Q_n/\prod_{j=2}^n(1-2/p_j)$ (see (\ref{indip-prob})) for each interval $I_{n-1}$ up to $n\simeq 10^4$. Here $Q_n$ is the  relative frequency of twin primes inside $I_{n-1}$. The (asymptotic) mean value $\bar \Theta \simeq 0.79371$, computed after removing the first $2\cdot 10^3$ entries,  is also reported as the horizontal solid line: we see that it is very close to the estimation obtained above. Addressing the whole computation to cousin primes leads to  pretty much the same result (thus we do not report the related picture), with a mean value $\bar \Theta \simeq 0.79357$. These results are consistent with the Hardy-Littlewood conjecture that twin and cousin primes share the same asymptotic distribution.

\subsection{New conjectures}
From Theorem \ref{mainthm} and \eqref{Qnt_lim} we see that the probability that the interval $I_{n-1}$ contains at least one twin prime is 
\begin{equation}
\label{qn}
q_n=1-(1-Q_n)^{T_n} \sim 1-\left(1-\frac{C_2 }{\log^2 p_n}\right)^{T_n},
\end{equation}
and the same holds true for cousin primes.
How fast $q_n$ converges to $1$ will depend on the rate of convergence to zero of the quantity $(1-Q_n)^{T_n}$. Since $T_n$ is proportional to the prime gap which has a random behavior (see for example \cite{Tabuguia}), we expect an irregular decay to zero of the quantity $(1-Q_n)^{T_n}$. This is displayed in the right picture of Figure \ref{FigQn}, which reveals a fan-shaped pattern where the slope of each {\em slat}  depends on the prime gap $p_n-p_{n-1}$. For example, the small circles in the uppermost branch correspond to those values of $n$ such that $(p_{n-1},p_n)$ forms a twin pair and thus $T_n=2p_n-2$. Analogously, the branch immediately below refers to the case where $(p_{n-1},p_n)$ is a  pair of cousin primes, so that $T_n=4p_n-8$, and so on for the subsequent branches. Inserting the asymptotic estimation of primes
$$
p_n\sim x(\log(x)+\log(\log(x))-1)
$$
in \eqref{qn},  we can get accurate analytical approximations of the branches: the solid line in the picture shows such an approximation for the uppermost branch. We see that for $n\ge 10^3$, the probability $q_n$ is greater than $1-10^{-50}$. This analytic behavior, together with a direct check of the existence of twin primes inside each $I_{n-1}$ up to $n=5\cdot 10^3$, suggests the following stronger conjectures on the existence of infinitely many twin and cousin primes.
\begin{conj}\label{Conjecture-tw} 
For any integer $n\geq 3$, there exists at least one pair of twin primes lying in the interval $(p_{n-1}^2,p_n^2)$. 
\end{conj} 
\begin{conj}\label{Conjecture-co} 
For any integer $n\geq 3$, there exists at least one pair of cousin primes lying in the interval $(p_{n-1}^2,p_n^2)$. 
\end{conj}


\section{Conclusion}\label{Sec5}
\noindent In the present work we have explored some theoretical properties hidden behind the algorithm  introduced in \cite{Bufalos}. The  paper provides a probabilistic approach to show the evidence of the twin and cousin prime conjectures.  Given any two consecutive primes $p_{n-1}$ and $p_n$, we have shown that there exists almost one twin and cousin prime  pair in the interval $I_{n-1}=(p_{n-1}^2,p_{n}^2)$ with probability $1$, as $n \to +\infty$. Though probabilistic in nature, we think that the interpretation of the algorithm as a discrete-time system inside each interval $I_{n-1}$, whose dynamics is ruled by a superposition of coprime periods, is prone to a deterministic analysis. A further path of  investigation concerns the quite natural generalization of the presented results to the case of primes that differ of a given even number.  These aspects  will be the subject of a future research.

\section*{Acknowledgements}
\noindent The third author wishes to thank  ``{\em Fondo acquisto e manutenzione attrezzature per la ricerca}'' D.R. 3191/2022, University of Bari.


\begin{thebibliography}{99} 
	\addcontentsline{toc}{chapter}{Bibliografia} 

\bibitem{Baker} R. C. Baker and T. Freiberg, Limit points and long gaps between primes, {\em Quarterly Journal of Mathematics}, 67 (2016), 233--260.

\bibitem{Banks} W. D. Banks, T. Freiberg and C. L. Turnage-Butterbaugh, Consecutive primes in tuples, {\em Acta
Arithmetica}, 167 (2015), 261--266. 


\bibitem{Bombieri} E. Bombieri, Problems of the Millennium: The Riemann Hypothesis, {\em Clay Mathematics Institute} (2000).

\bibitem{Brun} V. Brun, \"Uber das Goldbachsche Gesetz und die Anzahl der Primzahlpaare, {\em Archiv for Mathematik og Naturvidenskab} 34 (1915), 3--19.

\bibitem{Bufalos} D. Bufalo, M. Bufalo G. Orlando and R. Tetta, A new algorithm to find prime numbers with less memory requirements, {\em Journal of Discrete Mathematical Sciences and Cryptography} (2023), in press.

\bibitem{Castillo} A. Castillo, C. Hall, R. J. Lemke Oliver, P. Pollack and L. Thompson, Bounded gaps between
primes in number fields and function fields, {\em Proceeding of the American Mathematical Society} 143 (2015), 2841--2856.

\bibitem{Chen} J. R. Chen,  On the representation of a larger even integer as the sum of a prime and the
product of at most two primes, {\em Scientia Sinica}, 16 (1973), 157--176.


\bibitem{Dunham} W. Dunham, A note on the origin of the twin prime conjecture, {\em Notices of the International Consortium of Chinese Mathematicians}, 1 (2013), 63--65.

\bibitem{Erdos} P. Erd\"os, On the difference of consecutive primes, {\em Quarterly Journal of Mathematics} 6 (1935), 124--128. 

\bibitem{Freiberg} T. Freiberg, Short intervals with a given number of primes, {\em Journal of Number Theory}, 163 (2016), 159--171.


\bibitem{Glaisher}
J. W. L. Glaisher, An enumeration of prime-pairs, {\em Messenger of Mathematics}, 8 (1978),  28--33.


\bibitem{Gold1} D. Goldston, Y. Motohashi, J. Pintz, and C. Y. Y{\i}ld{\i}r{\i}m, Small gaps between primes exist, {\em Japan Academy Proceedings Series A. Mathematical Sciences} 82 (2006), 61--65.

\bibitem{Gold2} D. A. Goldston, S. W. Graham, J. Pintz, C. Y. Y{\i}ld{\i}r{\i}m, Small gaps between primes or almost primes, {\em Transactions of the American Mathematical Society} 361 (2009), 5285--5330.


\bibitem{Guy} R. K. Guy, {\em Unsolved problems in number theory}, Problem Books in Mathematics, Springer-Verlag, New York (2004).

\bibitem{Maier} H. Maier and M. T. Rassias, Large gaps between consecutive prime numbers containing
perfect k-th powers of prime numbers, {\em Journal of Functional Analysis}, 272 (2017), 2659--2696.

\bibitem{Maynard} J. Maynard, Small gaps between primes, {\em Annals of Mathematics} 181 (2015), 383--413.

\bibitem{Maynard2} J. Maynard, The twin primes conjecture, {\em Japan Journal of Mathematics} 14 (2019), 175--206. (2019);


\bibitem{Park} Y. Park and H. Lee, On the several differences between primes, {\em Journal of Applied Mathematics and Computing} 1 (2003), 3--51.

\bibitem{Polignac} A. de Polignac, Recherches nouvelles sur les nombres premiers, {\em Comptes Rendus} 29 (1849), 397--401.

\bibitem{Poly1} D. H. J. Polymath, New equidistribution estimates of Zhang type, {\em Algebra and Number Theory} 9 (2014),  2067--2199.

\bibitem{Poly2} D. H. J. Polymath, The ``bounded gaps between primes'' Polymath project: A Retrospective Analysis, {\em Newsletter of the European Mathematical Society} 94 (2014), 13--23.
 

\bibitem{Tabuguia} B. T. Tabuguia, An algorithmic random-integer generator based on the distribution of prime numbers, {\em Research Journal of Mathematics and Computer Science} 3 (2019).

\bibitem{Tao2011} T. Tao, Structure and randomness in the prime numbers, in: {\em An Invitation to Mathematics}
D. Schleicher, M. Lackmann (eds.), Springer-Verlag Berlin Heidelberg (2011).

\bibitem{Tao1} K. Ford, B. Green, S. Konyagin and T. Tao, Large gaps between consecutive prime numbers, {\em Annals of Mathematics}, 183 (2016), 935--974.

\bibitem{Tao2} K. Ford, B. Green, S. Konyagin, J. Maynard and T. Tao, Long gaps between primes, {\em Journal of
American Mathematical Society} 31 (2018), 65--105.

\bibitem{Tenenbaum} G. Tenenbaum,  Introduction to Analytic and Probabilistic Number Theory, Third Edition, Graduate Studies in Mathematics 163, American Mathematical Society Providence, Rhode Island (2015).

\bibitem{Vatwani} A. Vatwani, Bounded gaps between Gaussian primes, {\em Journal of Number Theory} 171 (2017), 449--
473.

\bibitem{Zhang} Y. Zhang, Bounded gaps between primes, {\em Annals of Mathematics}  179 (2014), 1121--1174.

\end{thebibliography}
\end{document}